\documentclass[final,a4paper,]{amsart}

\usepackage[]{hyperref}
\hypersetup{
pdftitle={High order semi-Lagrangian methods for the incompressible Navier--Stokes equations},
pdfauthor={Elena Celledoni, Bawfeh Kingsley Kometa, Olivier Verdier}}
\usepackage{caption}
\usepackage{subcaption}

\usepackage[usenames]{xcolor}
\usepackage{graphicx}
\usepackage{amssymb}
\usepackage{amsmath}
\usepackage{amsthm} 
\usepackage{mathrsfs}
\usepackage{mathtools}
\usepackage{bbm} 
\usepackage{multirow}

\usepackage{enumitem}

\newcommand*\loc[1]{\bar{#1}}
\newcommand*\glb[1]{#1} 

\usepackage{algorithm}
\usepackage{algpseudocode}
\usepackage{setspace}
\newenvironment{pseudocode}%
{\vspace{0.3cm}\begin{spacing}{1.5} \begin{algorithmic}}%
{\end{algorithmic} \end{spacing}\vspace{0.3cm}}

\newtheorem{remark}{Remark}[]

\newcommand{\ub}{\mathbf{u}}

\newcommand{\Vb}{\mathbf{V}}
\newcommand{\norm}[1]{\| #1\|}
\newcommand{\Xb}{\mathcal{X}}

\newcommand{\Nabla}{\boldsymbol{\nabla}}

\newcommand{\vf}[1]{\mathbf{#1}}
\renewcommand{\O}{\Omega}
\newcommand\T{\rule{0pt}{2.6ex}}
\newcommand\B{\rule[-1.2ex]{0pt}{0pt}}

\newcommand\inv{^{-1}}

\makeatletter
\long\def\@makecaption#1#2{%
  \vskip\abovecaptionskip
  \sbox\@tempboxa{#1: #2}%
  \ifdim \wd\@tempboxa >\hsize
   {\small\textbf{#1}: #2\par} 
  \else
    \global \@minipagefalse
    \hb@xt@\hsize{\box\@tempboxa\hfil}%
  \fi
  \vskip\belowcaptionskip}
\makeatother 
\makeatletter
     \def\fnum@table{\textbf{\tablename}\nobreakspace\textbf{\thetable}}
\makeatother

\usepackage{xparse}
\NewDocumentCommand\subfigc{ommm}{
\begin{subfigure}[rm]{#2}
                \centering
                \includegraphics[width=\textwidth]{#3}
                \caption{\small #4}\IfValueTF{#1}{\label{#1}}{}
\end{subfigure}}

\graphicspath{{./NSpaperFigsRevised/}}

\title[High order semi-Lagrangian methods]{High order semi-Lagrangian methods for the incompressible Navier--Stokes equations
}
\author{Elena Celledoni}
\author{Bawfeh Kingsley Kometa}
\author{Olivier Verdier}
\address{Institutt for matematiske fag, NTNU, 7049 Trondheim}

\begin{document}

\maketitle

\begin{abstract} 
We propose a class of semi-Lagrangian methods of high approximation order in space and time, based on spectral element space discretizations and exponential integrators of Runge--Kutta type. The methods were presented in \cite{ElenaBK2008} for simpler convection-diffusion equations. We discuss the extension of these methods to the Navier--Stokes equations, and their implementation using projections. Semi-Lagrangian methods up to order three are implemented and tested on various examples. The good performance of the methods for convection-dominated problems is demonstrated with numerical experiments.
\end{abstract}

\begin{description}
  \item[Mathematics Subject Classification (2010)]{Subject classification: Primary 54C40, 14E20; Secondary 46E25, 20C20.} 
\item[Keywords]
Navier--Stokes -- Projection -- Semi-Lagrangian -- Runge--Kutta
\end{description}

\section{Introduction}

Consider the incompressible Navier--Stokes equations
\begin{align}
\ub_t+\ub\cdot \Nabla \ub &=\nu \Nabla^2\ub-\nabla p,\label{NSvel}\\
\Nabla \cdot \ub &= 0.\label{NSinc}
\end{align}

\noindent Here $\ub=\ub(\mathbf{x},t)$ is the velocity field defined on the cylinder $\Omega \times [0,T]$ ($\Omega \subset \mathbf{R}^d$ for $d=2,3$), subject to the incompressibility constraint \eqref{NSinc}, while $p=p(\mathbf{x},t)$ is the pressure and plays the role of a Lagrange multiplier, and $\nu$ is the kinematic viscosity of the fluid. We consider no slip, or periodic boundary conditions when the domains allow it.

For no slip boundary conditions we will mostly consider the case 
\begin{equation}
\label{nosliphom}
\ub|_{\partial\Omega}=0.
\end{equation}
The variables $(\ub, p)$ are sometimes called \emph{primitive variables} and the accurate approximation of both these variables is desirable in numerical simulations. 

A typical approach for solving numerically convection-diffusion problems (the incompressible Navier--Stokes equations included) is to treat convection and diffusion separately, the diffusion with an implicit approach and the convection with an explicit integrator, see for example \cite{Canuto1988,Ascher1995,Ascher1997,KenCarp2003}. We will refer to these methods as \textit{implicit-explicit} methods (IMEX). 
The advantage of this approach is that  most of the spatial discretizations of the diffusion operator give rise to finite dimensional counterparts which are symmetric and positive definite, so the implicit integration of the diffusion requires only the solution of symmetric positive definite linear algebraic systems.

In this paper we propose high order discretization methods in time of semi-Lagrangian type, to be used in combination with high order spatial discretizations of the Navier--Stokes equations, as for example spectral element methods.
High order methods are particularly interesting when highly accurate numerical approximations of a given flow are required. An interesting field of application is the direct numerical simulation of turbulence phenomena, as pointed out for example in \cite{Karniadakis2001}. Another relevant situation is in connection with discontinuous-Galerkin methods, as an alternative to the use of explicit Runge--Kutta schemes. In this context, the purpose is to alleviate severe time-step restrictions imposed by CFL conditions. See for example \cite{Kanevsky2007} for the use of IMEX time-stepping schemes combined with discontinuous-Galerkin space discretizations, and \cite{Restelli2006} for semi-Lagrangian discontinuous Galerkin methods.

The integration methods we propose in this paper are implicit-explicit exponential integrators of Runge--Kutta  type. They combine the use of a \textit{diagonally implicit Runge--Kutta} method (DIRK) for the diffusion with a \textit{commutator-free exponential integrator}  (CF) for the nonlinear convection, and were denoted DIRK-CF in \cite{Elena2005,ElenaBK2008}. 
To achieve higher order, the nonlinear convection term needs to be approximated by a {\it composition} of linearised convection flows along constant convecting vector fields. Compared to the simpler convection diffusion problems treated in \cite{ElenaBK2008}, we here deal with the non-trivial, additional difficulty of enforcing the incompressibility constraint without compromising the order of the methods. This difficulty is present also when applying IMEX methods to the same problems, and our strategy to ensure incompressibility is equally valid for IMEX methods.

A significant advantage of the proposed schemes, compared for example to IMEX schemes, is that, because of the presence of exponentials of the linearised convection operator, they are amenable to semi-Lagrangian implementations. Similarly to other semi-Lagrangian methods (see e.g. \cite{Karniadakis2001}), the methods we proposed here allow for the use of considerably larger time-steps compared to Eulerian schemes, especially in convection dominated problems (i.e. at high Reynolds numbers in the Navier--Stokes equations). 

\subsection{Semi-Lagrangian features of the proposed exponential integrators}

To explain how the presence of exponentials of the convection operator can allow for semi-Lagrangian implementations, let us consider the simple linear convection 
model problem 
\begin{equation}
\label{simpleconvection}
\frac{Du}{Dt}
=0, 
\end{equation} 
where $\frac{Du}{Dt} \coloneqq \frac{\partial u}{\partial t} + \Vb\cdot\nabla $ is the total derivative of $u$,
with  $\mathbf{V}:\mathbf{R}^d\rightarrow \mathbf{R}^d$  the convecting vector field (which we assume to be not time dependent).  Approximations of simple convection problems of the type \eqref{simpleconvection} are highly relevant to the  approach we propose in this paper because they appear as building blocks in our methods. See section~\ref{section:methods}.

Let $h$ be a fixed step-size in time.
A simple example of a semi-Lagrangian scheme for \eqref{simpleconvection} is
\begin{displaymath}
\frac{u_{n+1}-u_n(\Xb(t_{n}))}{h}= 0,
\end{displaymath}
where $\Xb(t)$ is the characteristic path, solution (at time $t_n$) of the ordinary differential equation
\[
\dot{\Xb}(t)=\mathbf{V}(\Xb(t)), \qquad \Xb(t_{n+1})=\mathbf{x}, \qquad \mathbf{x}\in \Omega\subset \mathbf{R}^d
,
\]
and $t_{n+1}=t_n+h$.
The practical realization of this method requires:
\begin{itemize}
\item introducing a space discretization, where
$U_{n}$ is the numerical solution at time $t_n$ and on all nodes of the discretization grid $\Gamma$ (or $U_n$ belonging to a suitable finite element space);
\item an operator $\mathcal{I}({U_n})(\tilde{\mathbf{x}})$ interpolating $U_n$ and evaluating the result  on $\tilde{\mathbf{x}}\in \Omega$ (notice that $\tilde{\mathbf{x}}$ is not necessarily on the grid $\Gamma$); 
\item a suitable integration method to solve the equation
$\dot{\Xb}=\mathbf{V}(\Xb)$ backward in time and compute the characteristic paths; we denote by $\Phi_h^V(\mathbf{x})$ its numerical flow at time $t_n$ and with initial value $\mathbf{x}\in \Gamma$.
\end{itemize}
So the fully discrete method can be expressed in the form
\begin{equation}
\label{SLdiscrete1}
U_{n+1}=\mathcal{I}({U_n})({\Phi_h^V (\Gamma) }).
\end{equation}

We interpret $\mathcal{I}({U})({\Phi_h^V (\Gamma) })$ as an evolution operator in the sense of \cite{morton1984ggm}, \cite{childs1990cgm}. 
In fact, in an Eulerian perspective, a space semidiscretization of the convection 
problem \eqref{simpleconvection} 
yields
 $C(\cdot)$ a discrete convection operator and a system of linear ordinary differential equations (ODEs) of the type
\begin{equation}
\label{semidiscrete}
\dot{U}(t)=C(V)\,U(t), 
\end{equation}
where $V$ is a discrete version of $\mathbf{V}$, typically known only on the nodes of the discretization grid $\Gamma$. Assuming $U_n$ is the initial condition, the solution of \eqref{semidiscrete} at $t_{n+1}=t_n+h$ is
\begin{displaymath}
U_{n+1}=\exp(hC(V))\,U_n, 
\end{displaymath}
where $\exp$ denotes the exponential operator.
Unless stated otherwise in our methods we choose  $C(V)$  such that  
\begin{equation}
\label{interpexp}
\exp (hC(V))\, U_n=\mathcal{I}({U_n})({\Phi_h^V (\Gamma) })
\end{equation}
is satisfied\footnote{We also assume that the characteristic paths are  integrated to high accuracy.}.
This choice allows to view the semi-Lagrangian discretization as coming from a semi-discrete operator $C(V)$, and it simplifies the presentation of the algorithms in sections~\ref{DIRKCFprojNS:sec} and \ref{projDIRKCFNS:sec} where exponentials of the same type of \eqref{interpexp} enter as building blocks of the proposed methods.

\subsection{Error estimates of semi-Lagrangian methods allowing larger time steps.}

In the case of pure convection problems a well known error estimate for semi-Lagrangian methods, due to Falcone and Ferretti \cite{FerrettiFalcone1998,FerrettiFalconePrep}, gives a bound for the local error  $\tau(x_i,t_{n+1})$ of the type
\begin{displaymath}
|\tau(x_i,t_{n+1})|\le K\left( h^r+\frac{\Delta x^{q+1}}{h}\right),
\end{displaymath}
where $x_i$ is a generic grid-point and $t_{n+1}$ is time.
In this estimate, the term $h^r$ is the error due to the numerical approximation of the characteristics paths, the  term $\frac{\Delta x^{q+1}}{h}$ arises from the accumulation of the interpolation error, $r$ and $q+1$ are integers denoting the order of the time and space discretizations respectively, and $K$ is a constant independent on $h$ and $\Delta x$. This estimate of the error suggests that the spatial error is affected positively by the use of large time-steps $h$, moreover,  when high order interpolation is used (like in the case of high order spectral element methods), the integration of the characteristic paths should also be done at high accuracy and in particular an optimal $h$ could be chosen so that
\begin{displaymath}
h^r =\frac{\Delta x^{q+1}}{h}.
\end{displaymath}
These error estimates motivate the interest in designing semi-Lagrangian, exponential integrators achieving high order in time, when the adopted space-discretization is of high order.
In \cite{ElenaBK2008} and \cite{ElenaBK2010} we considered high order space discretizations by spectral element methods for convection-diffusion problems and provided high order semi-Lagrangian time-discretizations for nonlinear convection-diffusion problems. We also showed numerically that the proposed integrators do overcome nominal CFL stability restrictions.

So far the case of linear and nonlinear convection-diffusion equations have been considered. Semi-discretizations of the Navier--Stokes equations, giving rise to index 2 differential-algebraic systems, have been approached successfully by BDF-like multi-step  methods proposed in \cite{ElenaBK2010}. The connections of these methods to the  methods proposed in \cite{Maday1990}, \cite{Karniadakis2001} and \cite{Giraldo2003}, have also been explained.  Here we address the case of semi-Lagrangian methods based on one-step formulae, and more precisely Runge--Kutta type formulae.


%

\subsection{Projections and reformulation into ODE form to enforce incompressibility}
Given a time-stepping technique, a standard approach to adapt the method to the incompressible Navier--Stokes equations is by means of projections. 
The primary example of this technique, and most famous projection method for the incompressible Navier--Stokes equations is  the Chorin's projection method, proposed by Chorin \cite{Chorin1968,Chorin1969} and T\'emam \cite{Temam1969}.
%
The study of the temporal order of the Chorin's projection method was considered in \cite{Shen1992} and \cite{Rannacher1991} and it revealed order 1 in time for the velocity and only $\frac{1}{2}$ for the pressure. 
Such loss of order in time for the pressure is well known and 
has been analysed  for projection methods for Stokes and Navier--Stokes equations. Remedies to restore the full time order  are also known in the literature,  \cite{guermond06aoo}. 

In this paper we  choose the following  strategy. We first semi-discretize in space, taking care of boundary conditions, then project the equations at the space discrete level and eliminate the Lagrangian multipliers, and obtain a system of ordinary differential equations (ODEs).  Finally we
apply the exponential integrators to the resulting system of ordinary differential equations. The exponentials of the pure convection operators are approximated by computing characteristics and subsequent interpolation as in \eqref{interpexp}. 
We ensure incompressibility in two different ways described in \autoref{section:methods} and \ref{projDIRKCFNS:sec}. The numerical approximation of the pressure is obtained in a post-processing step. 

In  \autoref{section:methods} we present our main scheme and show how to obtain high order implicit-explicit and semi-Lagrangian methods for the 
incompressible Navier--Stokes equations. The formal order analysis of these integrators applied to ODE problems has been addressed in \cite{ElenaBK2008} and \cite{ElenaBK2009}, an extension of the order analysis and convergence of the methods to the PDE context, is outside the scope of the present paper.
In \autoref{section:Implementation} we describe the implementation details, we use techniques from \cite{Fischer1997} for the efficient solution of the linear algebraic systems. We obtain an overall strategy which resembles conventional projection schemes as described in \cite{guermond06aoo}, where,
at each step in time, one only needs to solve a sequence of decoupled elliptic equations for the velocity and the pressure.
In \autoref{NumExpr:sec} we report the numerical experiments. We provide numerical verification of the temporal order of the methods, and we demonstrate the clear benefits of the semi-Lagrangian approach in the case of convection-dominated problems.

\section{High order implicit-explicit and semi-Lagrangian methods of Runge--Kutta type for the incompressible Navier--Stokes equations}
\label{section:methods}

In this section we present the details of the high order integration schemes proposed in this paper.


After spatial discretization of \eqref{NSvel},  
we obtain a system of differential-algebraic equations of the type:
\begin{equation}
\begin{aligned}
B \dot{y}&= A\,y+C(y)\,y - D^Tz,\\
Dy&= 0.\\
\end{aligned}
\label{disNS}
\end{equation}

The matrices and vectors appearing in that equation have the following signification:
\begin{itemize}
   \item $A$ represents the discrete Laplacian;
   \item $C(y)$ is the discrete convection operator;
   \item $D$ is the discrete  divergence operator, so $D^T$ is a discrete gradient operator;
   \item $B$ is a mass matrix coming from the spatial discretization. We assume that $B$ is explicitly invertible, i.e., $B\inv$ is available;
   \item the vectors $y$ and $z$ represent numerical approximations of the velocity $\ub$ and the pressure $p$ respectively. 
\end{itemize}



\let\glborig\glb
\renewcommand\glb[1]{{#1}}

We can eliminate the Lagrangian multiplier $z$ from the system \eqref{disNS} by multiplying on both sides by $\glb{B}^{-1}$ and using the constraint $\glb{D}\dot{\glb{y}}=0$. This leads to the following system of ordinary differential equations (ODEs)
\begin{align}
\label{disNS4}
\glb{B} \dot{\glb{y}}&= \glb{A}\,\glb{y}+ \glb{C}(\glb{y})\,\glb{y} - \glb{B}\glb{H}\glb{B}^{-1}(\glb{A}+\glb{C}(\glb{y}))\glb{y},
\end{align}
where \begin{equation}\glb{H}\coloneqq\glb{B}^{-1}\glb{D}^T(\glb{D}\glb{B}^{-1}\glb{D}^T)^{-1}\glb{D}.\end{equation} 
We introduce the projection
\[
\glb{\Pi} \coloneqq I-\glb{H}
\]
allowing to write
the  ODE in the short form
\begin{equation}
\label{projeq}
 \dot{\glb{y}}=\glb{\Pi} \glb{B}^{-1}\glb{A}\,\glb{y}+ \glb{\Pi} \glb{B}^{-1}\glb{C}(\glb{y})\,\glb{y} .
 \end{equation}

In order to fully describe the methods, we need an IMEX method, as well as its DIRK-CF counterpart.
In this paper, we use the second and third order IMEX-RK schemes with stiffly-accurate and L-stable DIRK parts \cite{Ascher1997}.
We refer to them as IMEX2L and IMEX3L respectively.
We refer to the corresponding DIRK-CF methods as DIRK-CF2L and DIRK-CF3L respectively.
All the corresponding coefficients are given in \autoref{IMEX2Ltable} and \autoref{IMEX3Ltable}.

\begin{table}[htbp]\caption{IMEX2L and DIRK-CF2L coefficients }
\centering
$\gamma=(2-\sqrt{2})/2$ and $\delta=1-1/(2\gamma)$

$a_{i,j} \equiv\quad$
\begin{tabular}{r|ccc} \label{IMEX2Ltable} 
$0$  \\
$\gamma$ & $0$ & $\gamma$ \\
$1$ & $0$ & $ 1-\gamma$& $ \gamma $ \\\hline
	&$0$ & $ 1-\gamma $& $ \gamma $ \\ 
\end{tabular}
,
\qquad
$\hat{a}_{i,j} \equiv\quad$
\begin{tabular}{r|lll}
$0$  \\
$\gamma$ & $\gamma$ \\
$1$\B & $ \delta$ & $1-\delta $ & \\\hline
 & $\delta $\T & $1-\delta$ & $0$\\ 
\end{tabular}
\qquad $\alpha_{i,1}^k \coloneqq \hat{a}_{i,k}$
\end{table}

\begin{table}[htbp] \caption{IMEX3L and DIRK-CF3L coefficients}
\centering
$a_{i,j} \equiv\quad$
\begin{tabular}{r|ccccc} \label{IMEX3Ltable}
$0$  \\
$\frac{1}{2}$ & $0$ & $\frac{1}{2}$ \\
$\frac{2}{3}$\T & $0$ &  $\frac{1}{6}$ & $\frac{1}{2}$ \\
$\frac{1}{2}$\T & $0$ &  $-\frac{1}{2}$ & $\frac{1}{2}$ & $\frac{1}{2}$\\
$1$\T\B &  $0$ & $\frac{3}{2}$ & $-\frac{3}{2}$& $\frac{1}{2}$& $\frac{1}{2}$ \\\hline
 & $0$ &  $\frac{3}{2}$\T & $-\frac{3}{2}$ & $\frac{1}{2}$& $\frac{1}{2}$ \\ 
\end{tabular},
\qquad
$\hat{a}_{i,j} \equiv\quad$
\begin{tabular}{r|ccccc}
$0$  \\
$\frac{1}{2}$\T & $\frac{1}{2}$ \\
$\frac{2}{3}$\T  & $\frac{11}{18}$& $\frac{1}{18}$ \\
$\frac{1}{2}$\T & $\frac{5}{6}$&$-\frac{5}{6}$ & $\frac{1}{2}$\\
$1$\T\B& $\frac{1}{4}$& $\frac{7}{4}$& $\frac{3}{4}$& $-\frac{7}{4}$\\\hline
& $\frac{1}{4}$\T & $\frac{7}{4}$& $\frac{3}{4}$& $-\frac{7}{4}$ & $0$\\ 
\end{tabular}\\\vspace{2mm}
\quad
$\alpha_{i,1}^k = \hat{a}_{i,k} \qquad 1 \leq i, k \leq s,\quad s=5$

$\alpha_{s+1,1}^k = \begin{bmatrix} \frac{7}{12} & \frac{1}{5} & \frac{3}{4} & -\frac{2}{3} & 0 \end{bmatrix}$

$\alpha_{s+1,2}^k = \begin{bmatrix} -\frac{1}{3} & \frac{31}{20} & 0 & -\frac{13}{12} & 0\end{bmatrix}$
            
\end{table}

\subsection{DIRK-CF methods applied to the projected semi-discretized Navier--Stokes equations}\label{DIRKCFprojNS:sec}

In this section we present our main semi-Lagrangian approach to the incompressible Navier--Stokes equations based on the exponential integrators methods of \cite{ElenaBK2008}.
The operators $A$ and $D$  and $C$ are as in the previous section. 

We assume that we have chosen a DIRK-CF method (see \ref{section:IMEXDIRKCF}), so we have a Butcher tableau $a_{i,j}$ and coefficients $\alpha_{i,j}^k$.
We use the convention that $a_{s+1,s+1} = 0$.
In particular, we will use the coefficients in \autoref{IMEX2Ltable} and \autoref{IMEX3Ltable}.

Our main algorithm to solve \eqref{disNS} is \autoref{DIRKCFmethod}.

\begin{algorithm}
   \caption{DIRK-CF Method}
   \label{DIRKCFmethod}
   \begin{pseudocode}
   \For{$i \gets 1,s+1$}
   \State \begin{enumerate}
      \item  Compute the exponential  $\varphi_i$ (see section~\ref{sec:exponentials}): 
      \begin{align*}\varphi_i = \exp(h\glb{\Pi} \glb{B}^{-1}\glb{C}(\glb{Y}_i^{J}))\cdots\exp(h\glb{\Pi} \glb{B}^{-1}\glb{C}(\glb{Y}_i^{1}))\end{align*}
         where  \begin{equation*}\glb{Y}_i^{\gamma}\coloneqq\sum_{k=1}^{i-1} \alpha_{i,\gamma}^k\glb{Y}_k \quad\text{for}\quad\gamma=1,\dots , J\end{equation*}
      \item Solve in $Y_i, Z_i$ the linear equation:  \begin{align*}\glb{Y}_i - ha_{i,i}B\inv\glb{A}\glb{Y}_i &= \varphi_i\glb{y}_n + h\sum_{j=1}^{i-1}a_{i,j}\varphi_i\varphi_j^{-1}\glb{P_j}  + B\inv D^T Z_i\\
            D Y_i &= 0
         \end{align*}
      \item Solve in $P_i, Z_i'$ the linear equation: \begin{align*}P_i &= B\inv A Y_i + B\inv D^T Z'_i \\ D P_i &= 0 \end{align*}
   \end{enumerate}
     
   \EndFor
   \State $\glb{y}_{n+1}=\glb{Y}_{s+1}$
   \end{pseudocode}
\end{algorithm}

The rationale behind such an algorithm is that it is exactly equivalent to applying a standard DIRK-CF integrator to \eqref{projeq}:

\begin{pseudocode}
\For{$i=1:s+1$}
   \State  $\glb{Y}_i =\varphi_i\glb{y}_n+h\sum\limits_{j=1}^{i-1}a_{i,j}\varphi_i\varphi_j^{-1}P_j+ha_{i,i}\glb{\Pi} \glb{B}^{-1}\glb{A}\glb{Y}_i$
   \State $P_i = \glb{\Pi} \glb{B}^{-1}\glb{A}\glb{Y}_i$
   \State  $\glb{Y}_i^{\gamma}\coloneqq\sum_k \alpha_{i,\gamma}^k\glb{Y}_k$ for $\gamma=1,\dots , J$
 \State $\varphi_i = \exp(h\glb{\Pi} \glb{B}^{-1}\glb{C}(\glb{Y}_i^{J}))\cdots\exp(h\glb{\Pi} \glb{B}^{-1}\glb{C}(\glb{Y}_i^{1}))$ 
\EndFor
\State $\glb{y}_{n+1}=\glb{Y}_{s+1}$
\end{pseudocode}

\vskip0.2cm

\subsection{IMEX methods applied to the projected semi-discretized Navier--Stokes equations} \label{IMEXprojNS:sec}

The same strategy for enforcing the incompressibility constraint can be adopted for IMEX methods, we outline the details in this section.
IMEX methods will be used for comparison in the numerical experiments.


We need two Butcher tableaus $a$ and $\hat{a}$, for instance the ones in \autoref{IMEX2Ltable} or \autoref{IMEX3Ltable}.
We define the corresponding IMEX method applied to \eqref{projeq} in \autoref{IMEXmethod}.

\begin{algorithm}
   \caption{IMEX Method}
   \label{IMEXmethod}
   \begin{pseudocode}
   \For{$i \gets 1,s+1$}
   \State Solve for $Y_i, Z_i$ the linear equation 
   \begin{align*}
\glb{Y}_i - ha_{i,i}\glb{B}^{-1}\glb{A}\glb{Y}_i  &=   \glb{y}_n + h\glb{B}^{-1}\sum_{j=1}^{i-1}\left(a_{i,j}\glb{A}+\hat{a}_{i,j}\glb{C}(\glb{Y}_j)\right)\glb{Y}_j+ \glb{B}^{-1}\glb{D}^{T}Z_i\\
DY_i&=0
\end{align*}
   \EndFor
   \State $\glb{y}_{n+1}=\glb{Y}_{s+1}$
   \end{pseudocode}
\end{algorithm}

Note that this is equivalent to apply the IMEX Runge--Kutta method to the system \eqref{projeq}, as we would obtain.


\begin{pseudocode}
 \For{$i=1:s$}
   \State $\glb{Y}_i = \glb{y}_n + h\glb{\Pi}\glb{B}^{-1}\sum\limits_{j=1}^{i-1}\left(a_{i,j}\glb{A}+\hat{a}_{i,j}\glb{C}(\glb{Y}_j)\right)\glb{Y}_j + ha_{i,i}\glb{\Pi}\glb{B}^{-1}\glb{A}\glb{Y}_i$
\EndFor
\State $\glb{y}_{n+1} =\glb{y}_n + h\glb{\Pi}\glb{B}^{-1}\sum\limits_{i=1}^{s}\left(b_{i}\glb{A}+\hat{b}_{i}\glb{C}(\glb{Y}_i)\right)\glb{Y}_i $ 
%
\end{pseudocode}
 
See \ref{section:IMEXDIRKCF} for the definition of IMEX methods for convection-diffusion problems.


We notice in particular that 
$y_{n+1}$
satisfies the discrete incompressibility constraint $Dy_{n+1} = 0$ being the sum of terms which vanish when premultiplied by $D$.

Finally, to recover the correct approximation of the pressure at time $t_{n+1}$ we perform a post-processing step (which amounts to an extra projection). We consider the right hand side of \eqref{disNS} and evaluate it in $\glb{y}_{n+1}$ leading then to an approximation of  $\glb{D}\dot{\glb{y}}(t_{n+1})$. Since $\glb{D}\glb{y}_{n+1}=0$ and $\glb{D}\dot{\glb{y}}(t_{n+1})=0$ the correct approximation of the pressure is given by $z_{n+1}$ such that
\begin{displaymath}
\glb{B}^{-1}\glb{A}\glb{y}_{n+1}+\glb{B}^{-1}\glb{C}(\glb{y}_{n+1})\glb{y}_{n+1}-\glb{B}^{-1}\glb{D}^Tz_{n+1}=0.
\end{displaymath}
This amounts to solving a linear system for $z_{n+1},$ obtained by multiplying by $\glb{D}$:
\begin{equation}
\label{dPoiseq}
\glb{D}\glb{B}^{-1}\glb{D}^Tz_{n+1}=\glb{D}\glb{B}^{-1}(\glb{A}+\glb{C}(\glb{y}_{n+1}))\glb{y}_{n+1}.
\end{equation}


\section{Implementation issues}
\label{section:Implementation}

Before proceeding to the numerical experiments, we describe some of the implementation issues, related to the discretization of the Navier--Stokes equations, and to the use of spectral element methods. 

In the numerical experiments, the approximation is done in  $\mathbb{P}_{N}-\mathbb{P}_{N-2}$  compatible velocity-pressure discrete spaces. That is,  in each element we approximate the velocity by a $N$-degree Lagrange polynomial based on Gauss-Lobatto-Legendre (GLL) nodes in each spatial coordinate,  and the pressure by $(N-2)$-degree Lagrange polynomial based on  Gauss-Legendre (GL) nodes. The discrete spaces are  spanned by tensor product polynomial basis functions.
A consequence of
 this choice of the spatial discretization  is that the grid for the pressure does not include boundary nodes (there are no boundary conditions for the pressure). We remark that this is not the only viable choice of spatial discretization  for our time-integration schemes. 

\subsection{Computing the exponentials}
\label{sec:exponentials}

The exponential $\exp ( h \glb{\Pi} \glb{B}^{-1}\glb{C}(\glb{w}))\cdot g$ is the solution of
the semidiscretized equation
\begin{align}
\glb{B} \dot{\glb{v}}&=  \glb{C}(\glb{w})\,\glb{v} + \glb{D}^Tz,\\
\glb{D}\glb{v}&= 0, 
\end{align}
on the time interval $[0,h]$.

To approximate each of the exponentials $\exp ( h \glb{\Pi} \glb{B}^{-1}\glb{C}(\glb{w}))\cdot g$ we  use the following approach: 
we consider
\begin{equation}\label{SLproj_aprrox}
\exp ( h \glb{\Pi} \glb{B}^{-1}\glb{C}(\glb{w}))\cdot g=\glb{\Pi}\exp ( h \glb{B}^{-1}\glb{C}(\glb{w}))\cdot g+\mathcal{E}_h^q\cdot g+\mathcal{O}(h^{q+1}),
\end{equation}
where $\exp ( h  \glb{B}^{-1}\glb{C}(\glb{w}))\cdot g=\mathcal{I}(g)(\Phi_h^w(\Gamma))$ and 
\begin{equation}\label{SLproj_aprrox_correction}
\mathcal{E}_h^q := \sum_{k=2}^q\frac{h^k}{k!}((\glb{\Pi}B^{-1}C)^k-\glb{\Pi}(B^{-1}C)^k).
\end{equation}

For methods of order up to 2, it suffices to use the approximation
\begin{equation}\label{SLproj_aprrox_order2}
\exp ( h \glb{\Pi} \glb{B}^{-1}\glb{C}(\glb{w}))\cdot g= \glb{\Pi}\, \mathcal{I}(g)(\Phi_h^w(\Gamma))+\mathcal{O}(h^2).
\end{equation}
However for higher order methods, in general, values of $q\geq 2$ will be required for the correction operator \eqref{SLproj_aprrox_correction}. 

\begin{remark}
We observe from numerical tests that the more terms we include in the correction operator $\mathcal{E}_h^q $ defined in  \eqref{SLproj_aprrox_correction} the greater the accuracy of the methods. 
Nevertheless, only a few terms are required to achieve a desired order of convergence. 
For example, we observed that DIRK-CF3 methods constructed from the Butcher tableaus in  \autoref{IMEX3table:tab}, showed up to third order of convergence, even when the exponentials are approximated simply as in \eqref{SLproj_aprrox_order2} (with no additional correction term).
\end{remark}

\begin{table}[htbp] \caption{IMEX3 and DIRK-CF3 coefficients.}
\label{IMEX3table:tab}
 $\gamma = (3+\sqrt{3}\,)/6$

$a_{i,j} \equiv\quad$
\begin{tabular}{r|lll}
$0$ & $0$ & \\
$\gamma$ & $0$  & $\gamma$ & \\
$1-\gamma$ & $0$ & $1-2\gamma$ & $\gamma$\\\hline
 & $0$ & $\frac{1}{2}$ & $\frac{1}{2}$\T \\
\end{tabular}
\qquad
$\hat{a}_{i,j} \equiv\quad$
\begin{tabular}{r|lll}
$0$ & \\
$\gamma$ & $\gamma$ & \\
$1-\gamma$ & $\gamma-1$ & $2(1-\gamma)$ & \\\hline
 & $0$ & $\frac{1}{2}$ & $\frac{1}{2}$\T \\
\end{tabular}\\\vspace{2mm}
\begin{align*}
\alpha_{i,1}^k &= \hat{a}_{i,k},\quad  i = 1,\ldots,s,  \\
\alpha_{s+1,1}^k &= \begin{bmatrix}  x_1, & x_2, & (3x_1+3x_2-6c_2x_2-1)/(6c_3-3)\end{bmatrix},  \\
\qquad \alpha_{s+1,2}^k &= \begin{bmatrix}  x_2,  & \frac{1}{2}-x_2, & (6c_3-6x_1-6x_2+12c_2x_2-1)/(12c_3-6) 
 \end{bmatrix},  \\[2mm]
 \text{where }\quad
 c_1 &= 0, \quad  c_2 = \gamma, \quad  c_3 = 1-\gamma,\quad x_1 = 1/2,\quad x_2 = 1/3.
\end{align*}
\end{table}

\subsection{Pressure-splitting scheme}\label{pressspl:subsubsec}
This scheme is used to obtain a cost efficient computation 
of solutions of discrete linear Stokes systems (see e.g.,\cite{Fischer1997}). 

The IMEX and DIRK-CF methods described in sections~\ref{IMEXprojNS:sec} and ~\ref{DIRKCFprojNS:sec},  give rise to linear Stokes systems of the form
\begin{equation}\label{Stokes_stage}\left\{
\begin{aligned}
&\glb{\mathcal{H}}\glb{Y}_i - \glb{D}^TZ_i &&= \glb{B}\glb{f}_i\\
&\glb{D}\glb{Y}_i &&= \glb{g}_i
\end{aligned}\right.
\end{equation}
at each stage $i,$ where   $\glb{\mathcal{H}} =  \frac{1}{a_{i,i}h}\glb{B} - \glb{A},$ while $\glb{f}_i,\,\glb{g}_i $ incorporate 
the vector fields at earlier stage values $Y_j$(for $j<i$) and the contributions at boundary nodes (at stage $i$).
 
For a method of order 1 or 2, the pressure-splitting scheme can be carried out in the following steps:\footnote{Extension to higher order methods is straightforward (See \cite{Kvarving2010} and references therein).}

\begin{enumerate}[label=\textbf{Step \arabic*.}]
\item $\glb{\mathcal{H}}\hat{Y}_i - \glb{D}^Tz_n = \glb{B}\glb{f}_i$
\item  $\glb{D}\glb{B}^{-1}\glb{D}^T\Delta{z}_i = -\frac{1}{a_{i,i}h}(\glb{D}\hat{Y}_i-\glb{g}_i) $
\item  $\glb{Y}_i = \hat{Y}_i + a_{i,i}h\,\glb{B}^{-1}\glb{D}^T\Delta{z}_i, \quad Z_i = z_n + \Delta{z}_i.$
\end{enumerate}
The first step is an explicit approximation of the stage value of the velocity using the initial pressure $z_n.$ This approximation is not divergence-free. Steps 2 and 3 are thus the projection steps which enforce the algebraic constrain and correct the velocity and pressure. Note that this approximation introduces a truncation error of order $\mathcal{O}(h^3),$ and is thus sufficient for methods of order up to 2 (see e.g.\cite{Fischer1997}).  Solving \eqref{Stokes_stage} directly would lead to solving linear equations with the operator $\glb{D}\glb{\mathcal{H}}^{-1}\glb{D}^T$  for the $Z_i.$ However, the cost of inverting $\glb{D}\glb{\mathcal{H}}^{-1}\glb{D}^T$ is much higher than for inverting $\glb{D}\glb{B}^{-1}\glb{D}^T$ in Step 2, since $\glb{B}$ is usually diagonal or tridiagonal and easier to invert than $\glb{\mathcal{H}}$ (which is usually less sparse). This explains the main advantage for using the pressure-splitting schemes in the numerical computations. We have exploited this advantage in the numerical experiments  presented in sections~\ref{cavity:subsec} and \ref{shearrollup:subsec}.

The matrix $\glb{\mathcal{H}}$ is a discrete Helmholtz operator and is symmetric positive-definite (SPD); the mass matrix $\glb{B}$ is diagonal and SPD, and thus easy to invert. We use conjugate gradient methods for $\glb{\mathcal{H}}^{-1},$ with $\glb{B}^{-1}$ taken as preconditioner. The entire system \eqref{Stokes_stage} forms a symmetric saddle system, which has a unique solution for $\glb{Y}_i$ provided $\glb{D}$ is of full rank. The choice of spatial discretization method guarantees this requirement.  The system 
can be solved  by a Schur-complement approach (or block LU-factorization) and the pressure-splitting scheme.

Finally, we remark that the use of the pressure splitting scheme with our methods leads to overall approaches which can be regarded as a conventional projection schemes in the sense of \cite{guermond06aoo}.

\section{Numerical experiments}\label{NumExpr:sec}

For the numerical experiments we shall employ a spectral element method (SEM) based on the standard Galerkin weak formulation  as detailed out in \cite{Fischer2002}. See section \ref{section:Implementation}. 
We use a rectangular domain consisting of $N_e = N_x\times N_y$  uniform rectangular elements. The resulting discrete system has the form \eqref{disNS} or \eqref{eqlarge} (see section~\ref{section:methods}).  The semi-Lagrangian schemes associated to all the DIRK-CF methods in this section are achieved by tracing characteristics and interpolating as in \cite{Giraldo2003}. The fourth order explicit RK method is used for approximating the paths of characteristic.

\subsection{Temporal order tests for the IMEX methods}\label{ordertest:IMEX}
We investigate numerically the temporal order of convergence of the IMEX methods (contructed from  \autoref{IMEX2Ltable} and \autoref{IMEX3Ltable}) following the algorithm described in section~\ref{IMEXprojNS:sec}. 

\begin{figure}[t!]
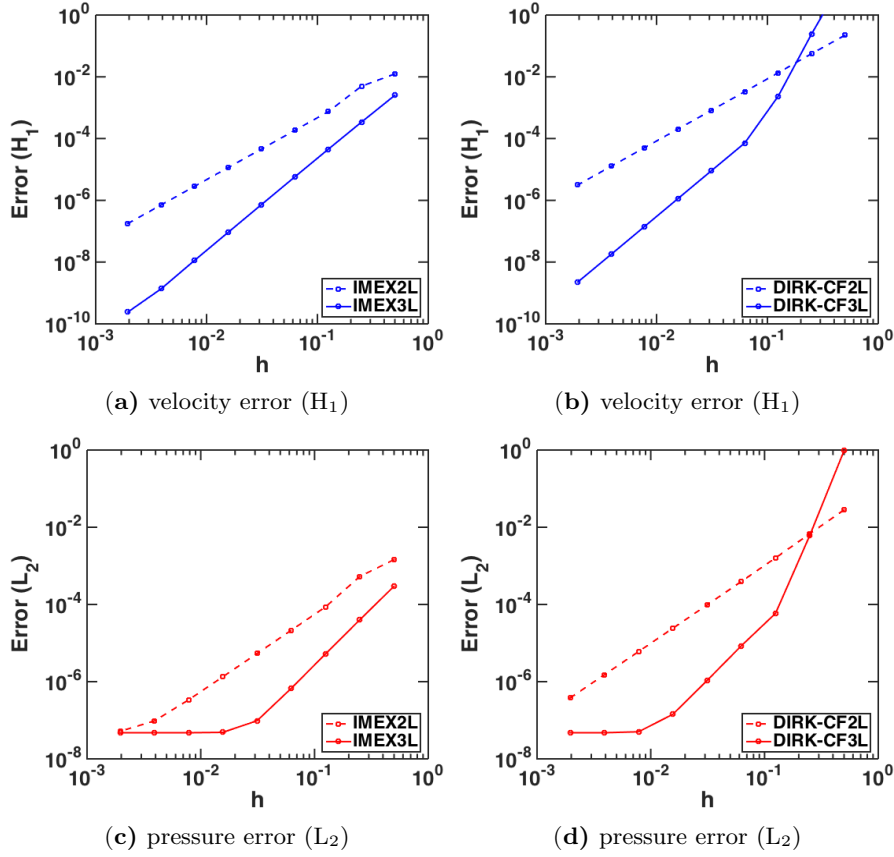

\centering
\subfigc[IMEXprojper:fig]{0.45\textwidth}{Order_temp_IMEX443_222projper}{velocity error (H$_1$)}
\subfigc[DIRKCFprojper2:fig]{0.45\textwidth}{Order_temp_DIRKCF443222projper2}{velocity error (H$_1$)}
\subfigc[IMEXprojperp:fig]{0.45\textwidth}{Order_temp_IMEX443_222projper_pp}{pressure error (L$_2$)}
\subfigc[DIRKCFprojper2p:fig]{0.45\textwidth}{Order_temp_DIRKCF443222projper2pp}{pressure error (L$_2$)}
\caption{Temporal convergence of IMEX2L, IMEX3L, DIRK-CF2L and DIRK-CF3L methods. Test problem: Taylor vortex (\ref{TaylorVortex}), with  $Re=2\pi^2$, $\O=[-1,1]^2$ and doubly-periodic bc. Discretization parameters:  $N=12$, $N_e=4$, $N_x = N_y = 2$, $h = T/2^k$, $k=1,\ldots,9$, $T=1$. }
\label{IMEXorderNS:fig}
\end{figure}

In the first example we consider the Taylor vortex problem  with exact solution  given by 
\begin{equation} \label{TaylorVortex}\left\{
\begin{aligned}
u_1(\mathbf{x},t) &= -\cos(\pi x_1)\sin(\pi x_2)\exp(-2\pi^2t/Re),\\
u_2(\mathbf{x},t)  &= \sin(\pi x_1)\cos(\pi x_2)\exp(-2\pi^2t/Re),\\
p(\mathbf{x},t)  &= -\frac{1}{4}[\cos(2\pi x_1) + \cos(2\pi x_2)]\exp(-4\pi^2t/Re),
\end{aligned}\right.
\end{equation}
where  $Re=1/\nu$ is the Reynolds number,  and $\ub\coloneqq (u_1,u_2),\,\,\mathbf{x}\coloneqq (x_1,x_2).$  The boundary condition is doubly-periodic on the domain $x_1,\,x_2\in[-1,1],$ and we choose $Re=2\pi^2.$ The initial conditions are determined from the exact solution \eqref{TaylorVortex}. For the spatial discretization we use a spectral method of high order $N=12,$ with $N_e = 4$ elements, and the time integration is done up to time $T=1.$ For each time-step $h = T/2^k,\,k=1,\ldots,9,$ the global error between the numerical solution and the exact PDE solution  (at time $T$) are measured in the $H_1-$ and $L_2$-norms\footnote{See \ref{Norms} for the definitions of these norms.} respectively, for the velocity and pressure. These are illustrated in log-log plots of the errors against the time-steps. The results for both the IMEX2L and IMEX3L show temporal convergence of order 2 and 3 respectively (see \autoref{IMEXorderNS:fig}). See also \autoref{IMEXorderNS:tab}.  
We notice that the measured global errors decrease as we decrease the time step.
However as the time step gets smaller,  the fixed spatial error becomes dominant  over the temporal error, and no further decrease in global error can be observed.
This is particularly the case for the pressure error of the third order method.
The spatial approximation space for the pressure is of lower order than the approximation space for the velocity (see the first paragraph of Section \ref{NumExpr:sec}).

\begin{table}[htbp]
\caption{Temporal convergence,  obtained with constant time steps $h = T/2^k$, where $k = 1,\ldots,9$, $T=1$. Test problem: Taylor vortex (\ref{TaylorVortex}).}
\label{IMEXorderNS:tab}
\centering
\begin{subfigure}[b]{\textwidth}
\centering
\caption{Approximation errors in the velocity, measured in the H$_1$-norm.}
\begin{tabular}{ccccc}
\hline 
 k & IMEX2L &  IMEX3L  & DIRK-CF2L & DIRK-CF3L \T \\ 
\hline 
1 & 1.2418e-02 &  2.5291e-03 &  2.2102e-01 &  1.8398e+01\T \\ 
2 & 4.9719e-03  &  3.4250e-04 & 5.5991e-02 &  2.4196e-01 \\ 
3 & 7.3930e-04 &  4.4798e-05 & 1.3504e-02 & 2.3406e-03 \\ 
4 & 1.8365e-04 & 5.7352e-06 & 3.3104e-03 & 7.2051e-05 \\ 
5 & 4.5775e-05 & 7.2573e-07 & 8.1987e-04 & 9.0946e-06 \\ 
6 & 1.1427e-05 & 9.1284e-08 & 2.0405e-04 & 1.1434e-06 \\ 
7 & 2.8545e-06 & 1.1445e-08 & 5.0903e-05 & 1.4337e-07 \\ 
8 & 7.1282e-07 & 1.3767e-09 & 1.2714e-05 & 1.7974e-08 \\ 
9 & 1.7843e-07 & 2.4201e-10 & 3.1763e-06 & 2.2809e-09 \\ 
\hline 
\end{tabular} 
~\\\vspace{0.2cm}
\end{subfigure}\\
\begin{subfigure}[b]{\textwidth}
\centering
\caption{Approximation errors in the pressure, measured in the L$_2$-norm.}
\begin{tabular}{ccccc}
\hline 
 k & IMEX2L &  IMEX3L  & DIRK-CF2L & DIRK-CF3L  \T \\ 
\hline 
1 & 1.4463e-03 &  2.9583e-04 &  2.8097e-02 &  9.5741e-01 \T \\ 
2 & 5.3052e-04  & 4.0100e-05 &  6.7856e-03 &  6.1495e-03 \\ 
3 & 8.6544e-05 &  5.2460e-06 &  1.6080e-03 &  5.8863e-05 \\ 
4 & 2.1504e-05 & 6.7327e-07 &  3.9240e-04 & 8.4854e-06 \\ 
5 & 5.3603e-06 & 9.7331e-08 &  9.7070e-05 & 1.0720e-06 \\ 
6 & 1.3389e-06 & 4.8632e-08 &  2.4152e-05 & 1.4275e-07 \\ 
7 & 3.3761e-07 & 4.7462e-08 & 6.0246e-06 & 5.0355e-08 \\ 
8 & 9.6013e-08 & 4.7443e-08 & 1.5054e-06 & 4.7489e-08 \\ 
9 & 5.1837e-08 & 4.7443e-08 & 3.7890e-07 & 4.7443e-08 \\ 
\hline 
\end{tabular} 
~\\\vspace{0.2cm}
\end{subfigure}\\
\begin{subfigure}[b]{\textwidth}
\centering
\caption{Convergence rates.\footnote{The approximate temporal convergence rates  are obtained by computing the linear slopes between any two adjacent time steps in the log-log plot of the global errors versus time steps.}}
\begin{tabular}{ccccccccc}
\hline 
 k & \multicolumn{2}{c}{IMEX2L} &  \multicolumn{2}{c}{IMEX3L}  & \multicolumn{2}{c}{DIRK-CF2L} & \multicolumn{2}{c}{DIRK-CF3L}  \T \\ 
   & vel. & press. & vel. & press. & vel. & press.  & vel. & press. \\ 
\hline 
1-2 & 1.3206 & 1.4469 & 2.8844 & 2.8831  & 1.9809 & 1.4469 & 6.2486 & 2.8831  \T \\ 
2-3 & 2.7496 & 2.6159 &  2.9346 & 2.9343 & 2.0518 & 2.6159 & 6.6918 & 2.9343\\ 
3-4 & 2.0092 & 2.0088 &  2.9655 & 2.9620 & 2.0283 & 2.0088 & 5.0217 & 2.9620\\ 
4-5 & 2.0043 & 2.0042 & 2.9823 & 2.7902 & 2.0136 & 2.0042 & 2.9859 & 2.7902\\ 
5-6 & 2.0021 & 2.0013 & 2.9910 & 1.0010 & 2.0064 & 2.0013 & 2.9917 & 1.0010\\ 
6-7 & 2.0011 & 1.9876 & 2.9957 & 0.0352 & 2.0031 & 1.9876 & 2.9955 & 0.0352\\ 
7-8 & 2.0016 & 1.8140 & 3.0554 & 0.0006 & 2.0013 & 1.8140 & 2.9958 & 0.0006\\ 
8-9 & 1.9982 & 0.8893 & 2.5081 & 0.0000 & 2.0010 & 0.8893 & 2.9782 & 0.0000\\ 
\hline 
\end{tabular} 
\end{subfigure}
\end{table}

%

\subsection{Temporal order tests for the  DIRK-CF methods}
\label{testDIRKCF:subsec}
Using the IMEX2L and IMEX3L methods, we construct two DIRK-CF methods, namely, DIRK-CF2L and DIRK-CF3L, of classical orders 2 and 3 respectively. 
Both DIRK-CF methods are applied to  \eqref{disNS} following the algorithm discussed in section~\ref{DIRKCFprojNS:sec}. To approximate the exponentials, we use a semi-Lagrangian approach coupled with a high order approximation based on \eqref{SLproj_aprrox}.
More precisely, we use  $\mathcal{E}_h^q$   (with $q=3$) for the DIRK-CF3L method,  and  \eqref{SLproj_aprrox_order2} for DIRK-CF2L.  In the case of DIRK-CF3L, larger time steps required higher values of $q>3$ in order to achieve convergence (that is, $q=3$ did not suffice).
The exponential on the right hand side of \eqref{SLproj_aprrox} is accurately computed using semi-Lagrangian methods for pure convection problems. 

Using the same test example as for section~\ref{ordertest:IMEX}, we observe the temporal order of convergence 2 and 3, in both the velocity and pressure (see \autoref{IMEXorderNS:fig}). See also  \autoref{IMEXorderNS:tab}.

\begin{figure}[htbp]
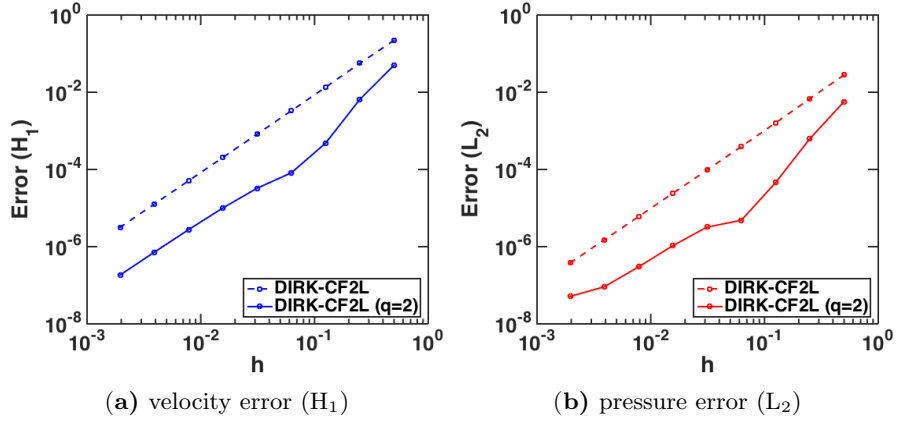

\centering
\subfigc[Order_temp_DIRKCF222projper2:fig]{0.45\textwidth}{Order_temp_DIRKCF222projper2}{velocity error (H$_1$)}
\subfigc[Order_temp_DIRKCF222projper2pp:fig]{0.45\textwidth}{Order_temp_DIRKCF222projper2pp}{pressure error (L$_2$)}
\caption{Temporal convergence of DIRK-CF2L method without correction (dashed line plot) and with correction    $\mathcal{E}_h^q$  (solid line plot). Test problem: Taylor vortex (\ref{TaylorVortex}), with  $Re=2\pi^2$, $\O=[-1,1]^2$ and doubly-periodic bc. Discretization parameters:  $N=12$, $N_e=4$, $N_x = N_y = 2$, $h = T/2^k$, $k=1,\ldots,9$, $T=1$. }
\label{DIRKCForderNS:fig}
\end{figure}

In addition, to show the impact of the correction operator   $\mathcal{E}_h^q$   on the accuracy of the DIRK-CF methods, we repeat the numerical test for DIRK-CF2L,  choosing this time $q=2$.  \autoref{DIRKCForderNS:fig} clearly show the improvement in accuracy.

\begin{figure}[t!]
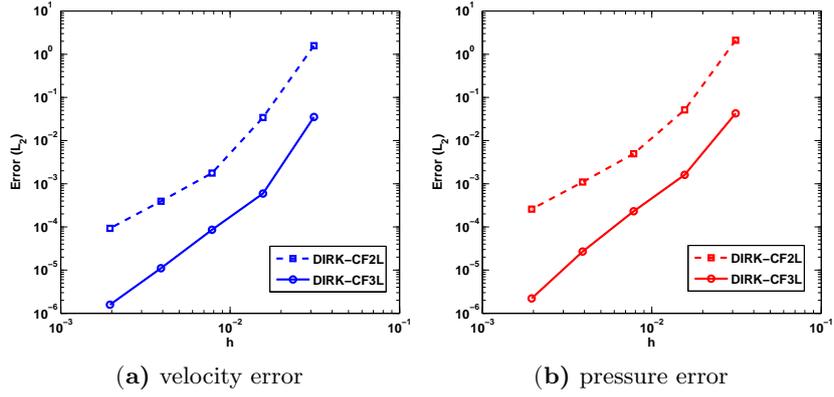

\centering
\subfigc{0.42\textwidth}{Order_temp_DIRKCF443222proj3}{velocity error}
\subfigc{0.42\textwidth}{Order_temp_DIRKCF443222proj3pp}{pressure error}
\caption{Temporal order of convergence. Test problem \eqref{Guermond2003}; $Re=100,\,T=1,\,N=8,\,N_e=16,\,N_x = N_y = 4,\,\O=[0,1]^2,\,h=h = T/2^k,\,k=5,\ldots,9.$ bc: homogeneous Dirichlet.  \textbf{(a)} velocity error ($L_2$):  DIRK-CF2L* (slope $=2.1268$), DIRK-CF3L* (slope $= 2.8632$);  \textbf{(b)} pressure error ($L_2$): DIRK-CF2L* (slope $=2.1296$),  DIRK-CF3L* (slope $=3.1628$).} 
\label{IMEXorderNS2:fig}
\end{figure}

To test the alternative semi-Lagrangian algorithm discussed in section~\ref{projDIRKCFNS:sec}, we apply the second and third order DIRK-CF methods on the test problem by \cite{GuermondShen2003} with exact solution given by 
\begin{equation} \label{Guermond2003}
\left\{
\begin{aligned}
u_1(\mathbf{x},t) &= \pi\sin(2\pi\,x_2)\sin^2(\pi\,x_1)\sin(t),\\
u_2(\mathbf{x},t) &= -\pi\sin(2\pi\,x_1)\sin^2(\pi\,x_2)\sin(t),\\
p(\mathbf{x},t) &= \cos(\pi\,x_1)\sin(\pi\,x_2)\sin(t),
\end{aligned}\right.
\end{equation}
for $x_1,\,x_2\in [0,1]$ and $t\in [0,T],$ with $T=1.$  A corresponding forcing term $\mathbf{f}$ is added to the momentum equations \eqref{NSvel} so that \eqref{Guermond2003} is the exact solution. In this test case we have used $Re=100.$ Meanwhile \eqref{Guermond2003} is used to prescribe the initial data and  boundary conditions (homogeneous Dirichlet on the entire boundary). The errors are all measured in the $L_2$-norm. We observe temporal order of convergence 2 and 3, in both the velocity and pressure (see \autoref{IMEXorderNS2:fig}). We rename the corresponding DIRK-CF methods by DIRK-CF2L* and DIRK-CF3L* respectively.

In the subsequent sections~\ref{cavity:subsec} and \ref{shearrollup:subsec}, we present a set of numerical experiments that illustrate the potential of the semi-Lagragian exponential integrators \cite{ElenaBK2008} for the treatment of convection-dominated problems. Two examples involving the incompressible Navier--Stokes models at high Reynolds numbers are considered. These examples  are the \textit{shear-layer roll up} problem in \cite{BrownMinion1995,FischerMullen2001,Fischer2002}, and the \textit{2D lid-driven cavity} problem (see \cite{Ghia2Shin1982,BrunSaad2006} and references therein). The second order  semi-Lagrangian DIRK-CF2L method (named SL2L in \cite{ElenaBK2008}) is used in each of these experiments.  The pressure-splitting  technique \cite{Fischer1997} (discussed in section~\ref{pressspl:subsubsec}) is applied to solve the discrete linear Stokes system that arises at each stage of the DIRK-CF method. 
The results reported in both sections~\ref{cavity:subsec} and \ref{shearrollup:subsec}  indicate that the semi-Lagrangian exponential integrators permit the use of large time-steps and Courant numbers. In both test cases the Courant number is defined following \cite[sect.5.2]{Giraldo1998} as
\begin{displaymath}
	C_r = \max\left(\frac{U\Delta t}{\Delta s}\right),
\end{displaymath}
where $U = \sqrt{\mathbf{u}_m^T\mathbf{u}_m}$ is the characteristic speed, and $\Delta s = \sqrt{\Delta{x}^2 + \Delta{y}^2}$ is the grid spacing. Here $\mathbf{u}_m$ denote the velocities at the midpoint of adjacent nodes.

\begin{figure}[b!]
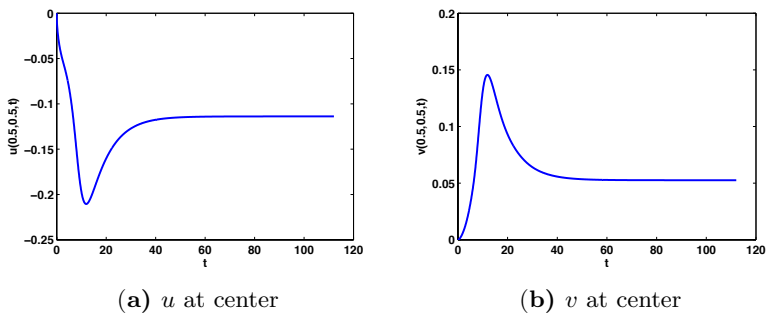

\centering
\subfigc{0.4\textwidth}{Cavity2DFlowRe400h03Ne1010p10_B}{$u$ at center}
\subfigc{0.4\textwidth}{Cavity2DFlowRe400h03Ne1010p10_B2}{$v$ at center}
\caption{Results of the second order semi-Lagrangian DIRK-CF method (SL2L) for the 2D lid-driven cavity problem. We have $(x,y)\in [0,1]^2$; $N_e=10\times 10$, $N = 10$, $h=0.03$, $Cr=9.0911$,  $Re=400$.
\textbf{(a)} Evolution of the horizontal velocity component $u$ at the domain center ($x=0.5$, $y=0.5$):  $t\in(0,112.08),$ \textbf{(b)} Evolution of the vertical velocity component $v$ at the domain center ($x=0.5, y=0.5$):  $t\in(0,112.08)$.}
\label{cavity2:fig}
\end{figure}


\subsection{Lid-driven cavity flow in 2D}\label{cavity:subsec}
We consider the 2D lid-driven cavity problem on a domain $(x,y)\in\Omega\coloneqq[0,1]^2$ with initial data $\vf{u} = (u,v) = (0,0)$ and constant Dirichlet boundary conditions
\begin{equation}
u = \begin{cases}
1 &\text{on upper portion of $\partial\Omega$}\\
0 &\text{elsewhere on $\partial\Omega$}
\end{cases}
,
\qquad v = 0 \quad\text{on}\quad\partial\Omega.
\end{equation}

\begin{figure}[t!]
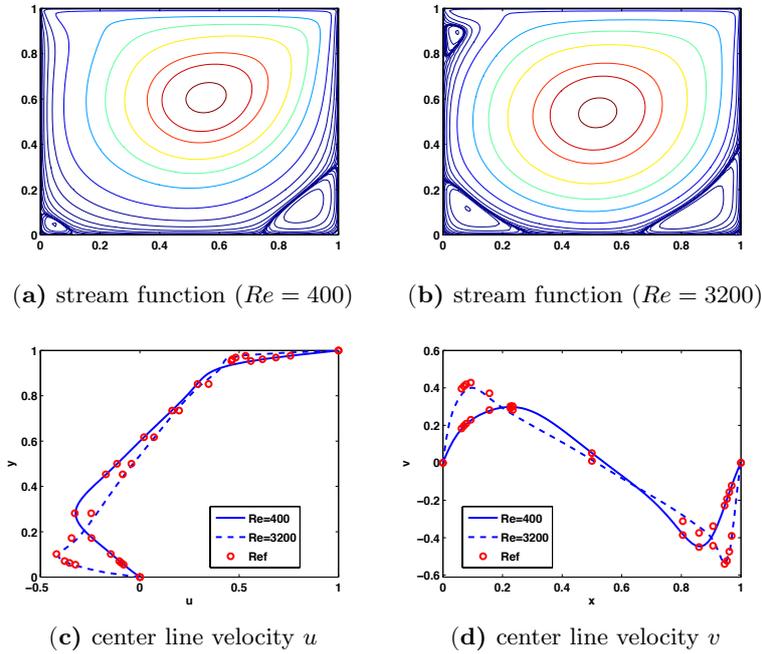

\centering
\subfigc[Cav2DRe400CF222_str1:fig]{0.4\textwidth}{Cav2DRe400CF222_streamlines1}{stream function ($Re=400$)}
\subfigc[Cav2DRe3200CF222_str1:fig]{0.4\textwidth}{Cav2DRe3200CF222_streamlines1}{stream function ($Re=3200$)}
\subfigc[Cav2DCF222_uc1:fig]{0.4\textwidth}{Cav2DCF222_uc1}{center line velocity $u$}
\subfigc[Cav2DCF222_vc1:fig]{0.4\textwidth}{Cav2DCF222_vc1}{center line velocity $v$}
\caption{Results of a second order DIRK-CF method for the 2D lid-driven cavity problem. We have $(x,y)\in [0,1]^2; \,N_e=10\times 10, \,N = 10, \,h=0.03, \,Cr=9.0911.$ In blue continuous line (our numerical solution); in red circles ({\Large \textcolor{red}{$\mathbf{\circ}$},} reference solution from \cite{Ghia2Shin1982}). \textbf{(a)} Streamline contours of the solution for $Re=400,$ \textbf{(b)} Streamline contours of the solution for $Re=3200$,
\textbf{(c)} Horizontal velocity component $u$ along the vertical center line ($x=0.5$), \textbf{(d)} Vertical velocity component $v$ along the horizontal center line ($y=0.5$). }
\label{cavity1:fig}
\end{figure}

We demonstrate the performance of the second order DIRK-CF method (SL2L, by the nomenclature of \cite{ElenaBK2008}). Spectral element method on a unit square domain $[0,1]^2$ with  $N_e=10\times 10$ uniform rectangular elements and polynomial degree $N=10$ is used (see \cite{Karniadakis2001}).  A constant time-step, $h=0.03,$ is used, corresponding to a Courant number of $Cr\approx 9.0911.$  
The time integration is carried out until the solution attains steady-state.  The  results in \autoref{cavity2:fig} show the evolution of the center velocity (at $Re=400$) up to steady state. It can be observed from this figure that steady state is attained at time $t\approx 40.$ At steady state the relative error ($L_2$-norm), between the velocity at a given time ($t_{n+1}$) relative to the velocity at the preceding time ($t_n$), has decreased to $\mathcal{O}(10^{-8}).$ 
The results also match with those of  \cite{Karniadakis2001}.
In \autoref{cavity1:fig}a-b we plot the streamline contours of the stream functions, choosing contour levels as in \cite{BrunSaad2006}. Meanwhile in \autoref{cavity1:fig}c-d plots of the  centerline velocities (continuous line, for $Re=400$, dashed line, for $Re=3200$) show a good match with those reported in \cite{Ghia2Shin1982} (plotted in red circles).

%

\subsection{Shear-layer roll up problem}\label{shearrollup:subsec}
We now consider the shear-layer problem \cite{BrownMinion1995,FischerMullen2001,Fischer2002} on a domain $\Omega\coloneqq[0,1]^2$ with initial data $\vf{u} = (u,v)$ given by
\begin{equation}
u = 
\begin{cases}
\tanh(\rho(y-0.25)) &\quad\text{for}\quad y\leq 0.5\\
\tanh(\rho(0.75-y)) &\quad\text{for}\quad y > 0.5
\end{cases}
,
\qquad v = 0.05\sin(2\pi\,x)
\end{equation}
which corresponds to a layer of thickness $\mathcal{O}(1/\rho).$
Doubly-periodic boundary conditions are applied. 

\begin{figure}[t!]
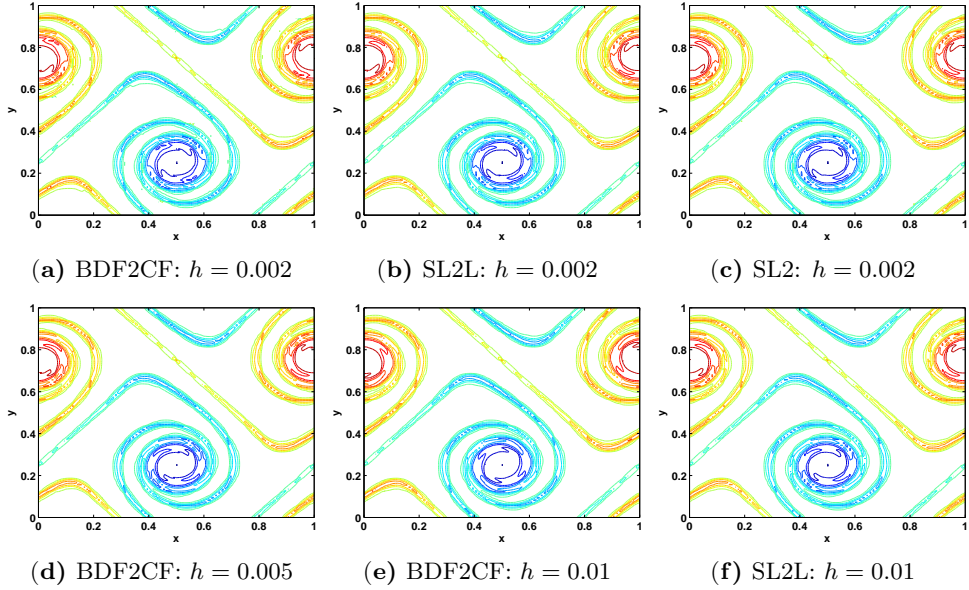

\centering
\subfigc{0.32\textwidth}{shearrollupBDF2CF002p8}{BDF2CF: $h=0.002$}
\subfigc{0.32\textwidth}{shearrollupCF222002p8}{SL2L: $h=0.002$}
\subfigc{0.32\textwidth}{shearrollupCF122002p8}{SL2: $h=0.002$}
\subfigc{0.32\textwidth}{shearrollupBDF2CF005p8}{BDF2CF: $h=0.005$}
\subfigc{0.32\textwidth}{shearrollupBDF2CF01p8}{BDF2CF: $h=0.01$}
\subfigc{0.32\textwidth}{shearrollupCF22201p8}{SL2L: $h=0.01$}
\caption{ Results of second order DIRK-CF methods (SL2 \& SL2L) and BDF2-CF method for the shear-layer rollup problem. We have $(x,y)\in [0,1]^2; \,N_e=16\times 16 = 256, \,N = 8.$ (filtering,  $ \alpha=0.3$),  $\rho=30, \,Re=10^5.$
Vorticity  contours (-70 to 70 by 15) of the solution at time $t=1.5.$ The corresponding Courant numbers are 
\textbf{(a)} {$Cr = 0.6393$}, 
\textbf{(b)} {$Cr = 0.6393$},
\textbf{(c)} {$Cr = 0.6393$},
\textbf{(d)} {$Cr = 1.5981$},
\textbf{(e)} {$Cr = 3.1963$},
\textbf{(f)} {$Cr = 3.1963$}.}
\label{shearrollup2:fig}
\end{figure}

In \autoref{shearrollup2:fig}  we demonstrate the performance of various second order methods including two DIRK-CF methods (SL2 \& SL2L, by the nomenclature of \cite{ElenaBK2008}), and also a second order semi-Lagrangian multistep exponential integrator (named BDF2-CF2, in \cite{ElenaBK2010}). The results are obtained at time $t=1.5,$  using a \emph{filter-based} spectral element method (see \cite{Fischer2002}) with $N_e=16\times 16$ elements and polynomial degree $N=8.$ The specified Reynolds number is $Re = 10^5,$ while $\rho=30$ and time-steps used are $h=0.002,\, 0.005,\,0.01$ corresponding to a Courant numbers of $Cr\approx 0.6393,\,1.5981,\,3.1963$ respectively. The filtering parameter used in each experiment is $\alpha=0.3$ (see for example \cite{Fischer2002}). However, the time-step and Courant number are up to about 10 times larger than that report in \cite{Fischer2002}. The initial values for the BDF2-CF are computed accurately using the second order DIRK-CF (SL2L) with smaller steps. The results are qualitatively comparable with those in \cite{FischerMullen2001,Fischer2002}.

%

\begin{figure}[t!]
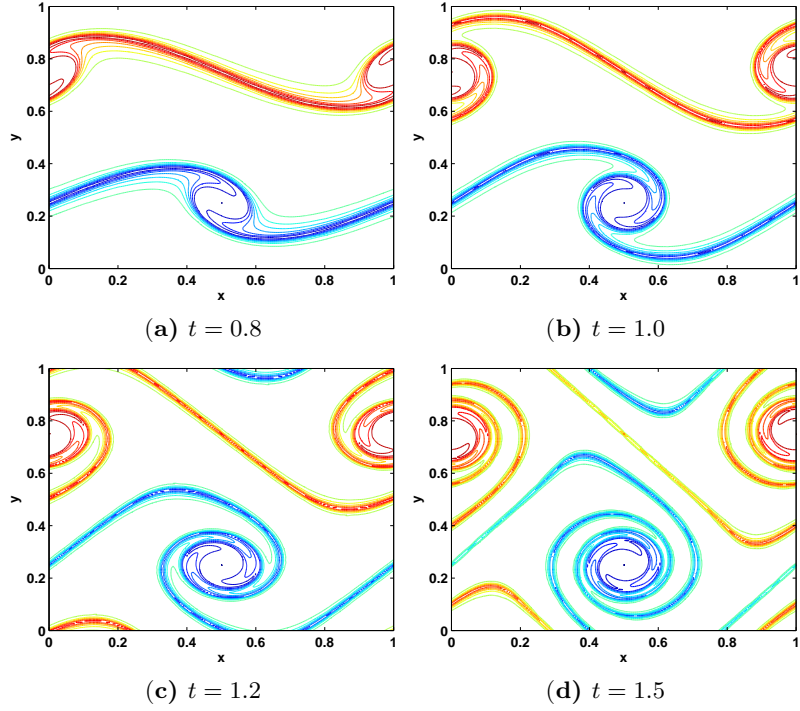

\centering
\subfigc{0.4\textwidth}{shearrollupCF222h01p16t08NoFilter}{$t=0.8$}
\subfigc{0.4\textwidth}{shearrollupCF222h01p16t1NoFilter}{$t=1.0$}
\subfigc{0.4\textwidth}{shearrollupCF222h01p16t12NoFilter}{$t=1.2$}
\subfigc{0.4\textwidth}{shearrollupCF222h01p16t15NoFilter}{$t=1.5$}
\caption{ Results of second order DIRK-CF method (SL2L) for the shear-layer rollup problem. We have $(x,y)\in [0,1]^2; \,N_e=16\times 16 = 256, \,N = 16, \,h=0.01,\,Cr = 11.9250,\,\rho=30, \,Re=10^5.$
Vorticity  contours (-70 to 70 by 15) of the solution at time \textbf{(a)} $t=0.8,$ \textbf{(b)} $t=1.0,$ \textbf{(c)} $t=1.2,$ \textbf{(d)} $t=1.5.$}
\label{shearrollup4:fig}
\end{figure}

In \autoref{shearrollup4:fig} we demonstrate the performance of the second order DIRK-CF method (SL2L). The results are obtained at times $t=0.8, 1.0, 1.2$ and $1.5$ respectively, using  spectral element method (without filtering) with $N_e=16\times 16$ elements and polynomial degree $N=16.$ The specified Reynolds number is $Re = 10^5,$ while $\rho=30.$ The time-step used is $h=0.01,$ corresponding to a Courant number of $Cr\approx 11.9250.$  This time-step is 10 times larger than that reported in \cite{Fischer2002}. Again the results are well comparable to those in  \cite{FischerMullen2001,Fischer2002}.

\begin{figure}[ht]
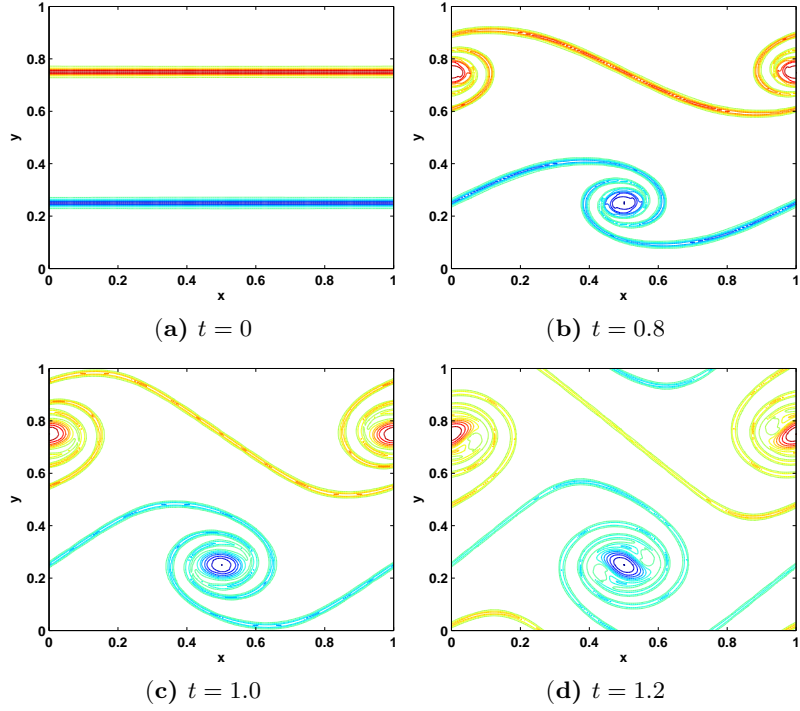

\centering
\subfigc{0.4\textwidth}{shearrollupthinCF222p16t0}{$t=0$}
\subfigc{0.4\textwidth}{shearrollupthinCF222h01p16t08NoFilter}{$t=0.8$}
\subfigc{0.4\textwidth}{shearrollupthinCF222h01p16t1NoFilter}{$t=1.0$}
\subfigc{0.4\textwidth}{shearrollupthinCF222h01p16t12NoFilter}{$t=1.2$}
\caption{Results of second order DIRK-CF method (SL2L) for the ``thin'' shear-layer rollup problem. We have $(x,y)\in [0,1]^2; \,N_e=16\times 16 = 256, \,N = 16, \,h=0.01,\,Cr = 11.9250.$ (\textbf{no filtering}),  $\rho=100, \,Re=40,000.$
Vorticity  contours (-36 to 36 by 13) of the solution at time \textbf{(a)} $t=0,$ \textbf{(b)} $t=0.8,$ \textbf{(c)} $t=1.0,$ \textbf{(d)} $t=1.2.$}
\label{shearrollup7:fig}
\end{figure}

Finally in \autoref{shearrollup7:fig} we demonstrate the performance of the second order DIRK-CF method (SL2L) for the ``thin'' shear-layer roll up problem, so defined for $\rho=100.$ The results are obtained at times $t=0.8, 1.0, 1.2$ and $1.5$ respectively, using  spectral element method (without filtering) with $N_e=16\times 16$ elements and polynomial degree $N=16.$ The specified Reynolds number is $Re = 4\times 10^ 4.$  The time-step used is $h=0.01,$ corresponding to a Courant number of $Cr\approx 11.9250.$  The results are well comparable to those in \cite{FischerMullen2001,Fischer2002}, except that we used 10 times the step size in time.

~

\section{Conclusion} In this paper we have presented a class of semi-Lagrangian methods for the incompressible Navier--Stokes equations of high order in time to be used with high order space discretizations, such as  for example spectral element methods. We have proposed a strategy to maintain the high temporal order also in the presence of constraints. As a by product, we have also derived projection methods based on IMEX Runge--Kutta schemes which have been used for comparison. The methods have been implemented and tested, and have been shown the predicted order of convergence in the case of periodic and no-slip boundary conditions. For convection-dominated test problems, in 2D with high Reynolds number, the semi-Lagrangian methods showed improved performance compared to their Eulerian counterparts, allowing for the use of considerably larger time-steps. So far methods up to order three have been implemented. Order four methods where obtained in \cite{ElenaBK2009}, and will be implemented and tested in future work.

\section*{Acknowledgements}
This work was supported in part by the GeNuIn project, grant from the Research Council of Norway.

%

\begin{appendix}
\section*{Appendix}
\setcounter{section}{1}

\subsection{IMEX and DIRK-CF methods for convection-diffusion equations}\label{section:IMEXDIRKCF}

A DIRK-CF method is a IMEX Runge--Kutta method of exponential type. It is defined by two Runge--Kutta tableaus one of implicit type to treat the linear diffusion and one of explicit type to treat the nonlinear convection by composition of exponentials.
These methods have a Runge--Kutta like format with two sets of parameters: 
$$\mathcal{A}=\{a_{i,j}\}_{i,j=1,\dots ,s},\quad \mathbf{b}=[b_1,\dots ,b_s],\quad \mathbf{c}=[c_1,\dots , c_s]$$ 
and 
$$\alpha_{i,l}^j,\quad \beta_l^{i},\quad i=1,\dots s, \,\, j=1,\dots ,s,\,\, l=1,\dots , J,\qquad \mathbf{\hat{c}}=[\hat{c}_1,\dots ,\hat{c}_s].$$
When applied to the convection-diffusion problem 
\begin{equation*}
\dot{U}(t)+C(U(t))\,U(t)=AU(t),
\end{equation*}
with linear diffusion and nonlinear convection the DIRK-CF methods have the following format:

\begin{pseudocode}
\For{i=1:s}
   \State $\mathcal{U}_i =\varphi_iU_n+h\sum\limits_{j=1}^{s}\varphi_{i,j}(a_{i,j}A\,\mathcal{U}_j)$
	 \State $\varphi_i = \exp(h\sum_k \alpha_{i,J}^kC(\mathcal{U}_k))\cdots\exp(h\sum_k\alpha_{i,1}^kC(\mathcal{U}_k))$ 
	 \State$\varphi_{i,j}\coloneqq\varphi_i \varphi_j^{-1}$
\EndFor
\State $U_{n+1}=\varphi_{s+1}U_n+h\sum_{i=1}^sb_{i}\varphi_{s+1,i}A\,\mathcal{U}_i$
\State $\varphi_{s+1}=\exp(h\sum_k\beta_J^kC(\mathcal{U}_k))\cdots\exp(h\sum_k\beta_{1}^kC(\mathcal{U}_k))$
\State $\varphi_{s+1,i}\coloneqq\varphi_{s+1}\varphi_i^{-1}$.
\end{pseudocode}

The methods are  associated to the two Butcher tableaus,
\begin{equation}
\label{tableaus}
\begin{array}{c|r}
   \mathbf{c}   &  \mathcal{A}\\[1mm]\hline
      & \mathbf{b}
\end{array},
\qquad
\begin{array}{c|r}
   \mathbf{\hat{c}}   &  \hat{\mathcal{A}}\\[1mm]\hline
      & \mathbf{\hat{b}}
\end{array},
\end{equation}
where we have defined
\begin{equation}
\label{eq:defcoeff}
\hat{a}_{i,j}\coloneqq\sum_{l=1}^J\alpha_{i,l}^j, 
\qquad
\hat{b}_{j}\coloneqq\sum_{l=1}^J\beta_{l}^j,
\end{equation}
for $i=1,\dots, s$, $\hat{\mathcal{A}}=\{\hat{a}_{i,j}\}_{i,j=1,\dots ,s}$ and $\mathbf{\hat{b}}=[\hat{b}_1,\dots ,\hat{b}_s]$.
The coefficients of the first tableau are used for the linear vector field $Ay$ while the coefficients of the second tableau, split up in the sums \eqref{eq:defcoeff}, are used for the nonlinear vector field $C(y)y$. 
We choose the first tableau to be a DIRK (diagonally implicit Runge--Kutta) method, this means we are solving only one linear system per stage.
The tableaus \eqref{tableaus}  are typically chosen so that they define a classical IMEX  method, which we call the underlying IMEX method, see \cite{Ascher1997} and \cite{ElenaBK2008} for more details. This IMEX method has the format

\begin{pseudocode}
\For{$i=1:s$}
   \State $\mathcal{U}_i =U_n+h\sum\limits_{j=1}^{s}(\hat{a}_{i,j}C(\mathcal{U}_j)\mathcal{U}_j+a_{i,j}A\,\mathcal{U}_j)$
\EndFor
\State $U_{n+1}=U_n+h\sum_{i=1}^s(\hat{b}_{i}C(\mathcal{U}_i)\mathcal{U}_i+b_{i}A\,\mathcal{U}_i).$
\end{pseudocode}
%

The order theory for classical IMEX methods reduces to the theory of partitioned Runge--Kutta methods, \cite{HLW2006}. Given an implicit and an explicit method of order $\kappa$ they must satisfy extra compatibility conditions in order for the corresponding IMEX method to have order $\kappa$. The extension of this theory to the DIRK-CF methods has been discussed in \cite{ElenaBK2008} and \cite{ElenaBK2009}.

\subsection{Definition of norms}\label{Norms}
For a square-integrable (respectively $H_1$) function $\mathbf{u} : \Omega \rightarrow \mathbb{R}^n,$ where $\Omega\subset\mathbf{R}^m$ is bounded and connected, the $L_2$-norm ($\norm{\cdot}_{L_2(\Omega)}$) and the  $H_1$-norm ($\norm{\cdot}_{H_1(\Omega)}$) are defined by
\begin{align}
\norm{\mathbf{u}}_{L_2(\Omega)} &\coloneqq  \left(\sum_{i=1}^n\int_{\Omega}u_i^2\,d\Omega\right)^{1/2},\\
\norm{\mathbf{u}}_{H_1(\Omega)} &\coloneqq  \left(\sum_{i=1}^n\int_{\Omega}(u_i^2+\nabla{u}_i\cdot\nabla{u}_i)\,d\Omega\right)^{1/2}.
\end{align}
In the spectral element approximations the continuous integrals of numerical solutions are accurately computed using Gauss quadrature rules.

\subsection{Boundary conditions and discrete stiffness summation}
\label{DSS}

For the sake of completeness, we illustrate the strategy for implementing the boundary conditions in the context of spectral element methods. We use the spectral element notion known as the direct-stiffness summation
(DSS), see for instance \cite{deville2002}.

Suppose 
we have to impose periodic or homogeneous Dirichlet boundary conditions, and that the variable $\loc{y}$ represents the values of the numerical solution at all discretization nodes in the computational domain (including boundary nodes). The variable $\glb{y}$ represents the restriction of $\loc{y}$ to the minimum degrees of freedom $k$ required to define the numerical solution, while $\loc{y}$ contains typically redundant components. Thus if the number of components of $\loc{y}$ is  $\mathcal{N},$ then $k<\mathcal{N}.$ We denote by $Q$ a prolongation or ``scatter'' operator such that $\loc{y} = Q\glb{y}.$ Associated to $Q$ is a restriction or ``gather'' operator denoted by $Q^T.$ The operator $Q$ is a $\mathcal{N}\times k$ constant matrix of rank $k\le \mathcal{N}$. 
The variable $\loc{y}$ is referred to as the \emph{local} variable, while $\glb{y}$ is the \emph{global} variable.
The DSS  operator $QQ^T$ ensures inter-element continuity and the fulfillment of the appropriate boundary conditions. So, for example, if the boundary conditions are periodic, given a vector $\loc{y}$ in the solution space or in the space of vector fields, $QQ^T\loc{y}$ is periodic.  


The spectral element discretization of the Navier--Stokes equations yields the discrete system  
\begin{align}
\label{eqlarge0}
\loc{B}\dot{\loc{y}}&= \loc{A}\loc{y}+\loc{C}(\loc{y})\loc{y}-\loc{D}^Tz\\
\loc{D}\loc{y}&= 0, \label{ceqlarge0}
\end{align}
where $\loc{y}$ is assumed to be in the range of $Q$ (i.e. $\loc{y} = Q\glb{y}$).
The relation between the local and global operators is $\glb{B}=Q^T\,\loc{B}\,Q$, $\glb{A}=Q^T\,\loc{A}\,Q$, $\glb{C}(\glb{y})=Q^T\, \loc{C}(Q\glb{y})\,Q$ and $\glb{D}=\loc{D}\,Q$.

Applying on both sides of \eqref{eqlarge0} the DSS operator $QQ^T$ we obtain
\begin{align}
\label{eqlarge}
\Sigma \dot{\loc{y}}&= QQ^T \loc{A}\loc{y}+QQ^T\loc{C}(\loc{y})\loc{y}-QQ^T\loc{D}^Tz,\\
\loc{D}\loc{y}&= 0.
\end{align}
The matrix $\Sigma=QQ^T\loc{B}$  is $\mathcal{N}\times \mathcal{N}$ and invertible on the range of $Q$. In practice the integration methods are reformulated for the local variable $\loc{y}$ and the local operators.

Indeed, in the computations the full data for the local variable $\loc{y}$ is stored, since all computations involving  the operators $\glb{B}^{-1}, \glb{\mathcal{H}}^{-1}$ and $(\glb{D}\glb{B}^{-1}\glb{D}^T)^{-1}$ must be done within the range of $Q.$ These operators are symmetric and positive-definite, and so they can be inverted using a fast iterative solver such as the conjugate gradient method. For example the problem $\glb{y} = \glb{\mathcal{H}}^{-1}\glb{f}$ is reformulated as follows:
Find $\loc{y}$ such that 
\begin{equation}
\left\{\begin{aligned}
QQ^T\mathcal{\loc{H}}\loc{y} &= Q\glb{f}\\
Q\glb{y} &= \loc{y}
\end{aligned}\right.,
\end{equation}
where $\mathcal{\loc{H}} = \frac{1}{a_{i,i}h}\loc{B} - \loc{A}.$

We refer to \cite{Fischer1997} for further details on DSS and boundary conditions. In the experiments reported in this paper, no special treatment has been taken to enforce pressure boundary conditions, since the discrete pressure space is not explicitly defined on discretization nodes on the boundary.



\subsection{Reformulation of the integration methods in local variables}

In this section we briefly discuss the correct implementation of the methods of \autoref{section:methods} as applied to \eqref{eqlarge}.
The purpose of reporting here these implementation details is to explicitly highlight when care has to be taken in the implementation. See also the remark below.

We first perform the elimination of the discrete pressure Lagrangian multiplier from  \eqref{eqlarge} in analogy to \eqref{disNS4}. 
We use 
$$\loc{H}=\loc{D}^T(\glb{D}\glb{B}^{-1}\glb{D})^{-1}\glb{D}\glb{B}^{-1}Q^T,$$ 
and we  get a system of ODEs for the variable $\loc{y}$:
\begin{align}
\label{disNS5}
\Sigma \dot{\loc{y}}&= QQ^T\loc{A}\,\loc{y}+QQ^T \loc{C}(\loc{y})\,\loc{y} - QQ^T\loc{H}(\loc{A}\,\loc{y}+\loc{C}(\loc{y})\,\loc{y}) ,
\end{align}

Introducing the projection
$\loc{\Pi}=I-\loc{H}$ allows us to write the following projected system of ODEs
\begin{equation}
\label{dNSproj}
\Sigma \dot{\loc{y}}=
QQ^T\loc{\Pi} \loc{A}\loc{y}+ QQ^T\loc{\Pi} \loc{C}(\loc{y})\loc{y}.
\end{equation}

\vskip0.3cm
Denote $\tilde{\Pi} = Q\glb{\Pi}\glb{B}^{-1}Q^T.$ Applying the method of section~\ref{IMEXprojNS:sec} to the ordinary differential equation \eqref{dNSproj} split in its projected convection and projected diffusion terms, and then rewriting it as a method for the differential-algebraic equation \eqref{eqlarge} we obtain

\begin{pseudocode}
\For{$i=1:s$}
   \State $\left(\Sigma - ha_{i,i}QQ^T\loc{A}\right)\loc{Y}_i + ha_{i,i}QQ^T\loc{D}^TZ_i = \Sigma\,\loc{y}_n + h\, \Sigma\,\tilde{\Pi}\sum\limits_{j=1}^{i-1}\left(a_{i,j}\loc{A} + \loc{C}(\loc{Y}_j)\right)\loc{Y}_j $
   \State $\loc{D}\loc{Y}_i = 0$
\EndFor
\State $\loc{y}_{n+1} = \loc{Y}_s,$
\end{pseudocode}

\noindent under the assumption $\loc{y}_n = Q\glb{y}_n.$ 
Since $\loc{y}_{n+1}=\loc{Y}_s,$ the approximation of the velocity satisfies the discrete incompressibility constraint, $\loc{D}\loc{y}_{n+1}=0$.

\begin{remark}\normalfont
We observe that this is not equivalent to what we obtain applying directly the IMEX method to \eqref{eqlarge}, which written in ODE form is
\begin{displaymath}
\Sigma \dot{\loc{y}}=QQ^T \loc{A}\loc{y}+QQ^T\loc{C}(\loc{y})\loc{y}-QQ^T\loc{H}(\loc{A}y+\loc{C}(\loc{y})\loc{y}).
\end{displaymath}
In fact if we apply the IMEX method to \eqref{eqlarge} we need to treat the term $QQ^T\loc{D}^Tz$ either with the implicit method or with the explicit method, while in the approach outlined in this section we treated $QQ^T\loc{H}\loc{C}(\loc{y})\loc{y}$ explicitly and $QQ^T\loc{H}\loc{A}\loc{y}$ implicitly, see also \cite{Kometa2011}. 
\end{remark}
\vskip0.3cm

Analogously, the method of section~\ref{DIRKCFprojNS:sec} applied to for \eqref{dNSproj} becomes 

\begin{pseudocode}
\For{$i=1:s+1$}
 \State {\footnotesize $ (\Sigma-ha_{i,i}QQ^T\loc{A})\loc{Y}_i+ha_{i,i}QQ^T\loc{D}^TZ_i=\Sigma Q\varphi_i\glb{y}_n+h\sum\limits_{j=1}^{i-1}a_{i,j}\Sigma\,Q\varphi_i\varphi_j^{-1}\,\loc{A}  \loc{Y}_j$}
   \State $\loc{D}\loc{Y}_i=0$ 
    \State  $\glb{Y}_i^{\gamma}\coloneqq\sum_k \alpha_{i,\gamma}^k\glb{Y}_k$ for $\gamma=1,\dots , J$
 \State $\varphi_i = \exp(h\glb{\Pi} \glb{B}^{-1}\glb{C}(\glb{Y}_i^{J}))\cdots\exp(h\glb{\Pi} \glb{B}^{-1}\glb{C}(\glb{Y}_i^{1}))$ 
\EndFor
\State $\loc{y}_{n+1}=\loc{Y}_{s+1}.$
\end{pseudocode}

The method of section~\ref{projDIRKCFNS:sec} may be reformulated in a similar way.
%

\subsection{Projected DIRK-CF methods for Navier--Stokes equations} \label{projDIRKCFNS:sec}

In this section we present an alternative semi-Lagrangian approach compared to the previous section~\ref{DIRKCFprojNS:sec}. Also in this case the integration methods  are a variant of the exponential integrators methods of \cite{ElenaBK2008}.

We consider \eqref{projeq}, 
and 
rearrange the terms in the form
 \begin{equation}
\label{projeqbroken1}
 \dot{\glb{y}}= \glb{B}^{-1}\glb{A}\,\glb{y}-\glb{H}\glb{B}^{-1}(\glb{A}+\glb{C}(\glb{y}))\,\glb{y}+  \glb{B}^{-1}\glb{C}(\glb{y})\,\glb{y}.
 \end{equation}

 This is a differential equation on the subspace of discrete divergence-free vector fields,  
 i.e. $\glb{D}\glb{y}=0$ for all $t$. To approximate the solution of this equation, we consider a projection method of the type reviewed in  \cite[IV.4]{HLW2006}, see also \cite[Sect.5.3.3]{SoellnerFuhrer1998}  and \cite[Sect.VII.2]{HW2ED}. The idea is to use a one-step integrator $\phi_h$ for advancing the numerical solution of \eqref{projeqbroken1} by one step, and an orthogonal projection on the subspace of divergence-free vector fields applying $\glb{\Pi}$ at the end of each step.
We choose $\phi_h$ to be the following integration method, in which  the coefficients of both the DIRK-CF method and the underlying IMEX method are used:
\begin{itemize}
\item the term $(I-\glb{H})\glb{B}^{-1}\glb{A}\,\glb{y}$ is treated implicitly with the DIRK coefficients, 
\item the term $\glb{H}\glb{B}^{-1}\glb{C}(\glb{y})\,\glb{y}$ is treated explicitly with the coefficients of the underlying explicit method,
\item  the term $\glb{B}^{-1}\glb{C}(\glb{y})\,\glb{y}$ is treated with the coefficients of the corresponding CF method.
\end{itemize}
The projection $\glb{\Pi}$ is used to guarantee divergence-free numerical approximations, i.e. $\glb{D}\glb{y}_{n+1}=0$. We obtain

\begin{pseudocode}
\For{$i=1:s+1$}
    \State$\glb{Y}_i =\varphi_i\glb{y}_n+h\sum\limits_{j=1}^{i-1}\varphi_i\varphi_j^{-1} \left(a_{i,j}\glb{\Pi}\glb{B}^{-1}\glb{A}\glb{Y}_j-\hat{a}_{i,j}\glb{H}\glb{B}^{-1}\glb{C}(\glb{Y}_j)\glb{Y}_j\right)+ha_{i,i}\glb{\Pi}\glb{B}^{-1}\glb{A}\glb{Y}_i$
   \State $\glb{Y}_i^{\gamma}\coloneqq\sum_k \alpha_{i,\gamma}^k\glb{Y}_k$ for $\gamma=1,\dots , J$
 \State$\varphi_i = \exp(h\glb{B}^{-1}\glb{C}(\glb{Y}_i^{J}))\cdots\exp(h\glb{B}^{-1}\glb{C}(\glb{Y}_i^{1}))$ 
\EndFor
\State $\glb{y}_{n+1}=\glb{\Pi}\glb{Y}_{s+1}$,
\end{pseudocode}

\noindent where 
$$a_{s+1,j} = b_j,\,\,\hat{a}_{s+1,j} = \hat{b}_j,\,\,\alpha_{s+1,\gamma}^k = \beta_{\gamma}^k,\,\,j,k = 1,\ldots,s,\,\, a_{s+1,s+1} = 0,\,\,\hat{a}_{s+1,s+1} = 0.$$


Since $\glb{\Pi}$ is an orthogonal projection, the order of the method is not affected by the use of the projection in the update step\footnote{The projection does not need to be orthogonal, but should be guaranteed not to  compromise the order of the method, an orthogonal projection will have this property.  The target of the orthogonal projection map on the discrete divergence-free subspace is the element of shortest distance to the point which is projected. Since $\glb{y}_{n+1}=\glb{\Pi}\glb{Y}_s$ and 
 $\|\glb{y}(t_{n+1})-\glb{Y}_s\|= \mathcal{O}(h^{r+1}),$ where $r$ is the order of the integration method, then 
 $$\|\glb{y}(t_{n+1})-\glb{y}_{n+1}\|\le\|\glb{y}(t_{n+1})-\glb{Y}_s\|+\|\glb{Y}_s-\glb{y}_{n+1}\|= \mathcal{O}(h^{r+1}),$$
because 
$$\|\glb{y}_{n+1}-\glb{Y}_s\|\le\|\glb{y}(t_{n+1})-\glb{Y}_s\|.$$}.
 This approach requires only exponentials of pure convection problems, this will ease the implementation of the method as a semi-Lagrangian method, because now we simply use the semi-Lagrangian approximation
 \begin{displaymath}
\exp\left(h\glb{B}^{-1} \glb{C}\big(\glb{w}\big)\right)\,g= \mathcal{I}(g)(\Phi_h^w(\Gamma)).
\end{displaymath}

Observe that at each stage $\glb{Y}_i $ does not necessarily satisfy $\glb{D}\glb{Y}_i =0$. 

\end{appendix}

\bibliographystyle{amsplain}
\bibliography{MyReviewedArticlesII}

\end{document}